\documentclass [11pt,oneside,mathscr]{amsart}
\usepackage{amscd,amsmath,amssymb,euscript}
\usepackage[frame,cmtip,curve,arrow,matrix,line,graph]{xy}

\title{Derived categories of coherent sheaves}
\author{Tom Bridgeland}

\date{}

\newtheorem{thm}{Theorem}[section]
\newtheorem{conj}[thm]{Conjecture}
\newtheorem{prob}[thm]{Problem}
\newtheorem{cor}[thm]{Corollary}
\newtheorem{prop}[thm]{Proposition}

\theoremstyle{definition}
\newtheorem{defn}[thm]{Definition}

\newtheorem{example}[thm]{Example}

\newcommand{\Pol}{\operatorname{\mathcal{P}}}
\newcommand{\NS}{\operatorname{NS}}

\renewcommand{\leq}{\leqslant}

\newcommand{\ch}{\operatorname{ch}}

\newcommand{\lEnd}{\operatorname{\mathcal{E}nd}_{\OO_Y}}
\newcommand{\lHom}{\operatorname{\mathcal{H}om}}
\newcommand{\Per}{\operatorname{Per}(Y/X)}
\newcommand{\Perw}{\operatorname{Per}(Y/Y')}

\newcommand{\D}{\operatorname{{D}}}
\newcommand{\Coh}{\operatorname{Coh}}
\newcommand{\Aut}{\operatorname{Aut}}
\newcommand{\isom}{\cong}

\newcommand{\tensor}{\otimes}
\newcommand{\PP}{\operatorname{\mathbb P}}
\newcommand{\C}{\mathbb C}

\newcommand{\RR}{\mathbb{R}}

\newcommand{\F}{\mathcal F}
\newcommand{\mat}[2]{\left( \begin{array}{cc} #1 \\ #2 
\end{array} \right)}
\newcommand{\rk}{\operatorname{rank}}
\newcommand{\dg}{\operatorname{deg}}
\newcommand{\Z}{\mathbb Z}
\newcommand{\A}{\mathcal A}

\newcommand{\OO}{\mathcal O}
\newcommand{\into}{\hookrightarrow}

\newcommand{\Pic}{\operatorname{Pic}}
\newcommand{\Hom}{\operatorname{Hom}}
\newcommand{\eu}{\operatorname{\chi}}

\newcommand{\lRa}[1]{\xrightarrow{\ #1\ }}
\newcommand{\lra}{\longrightarrow}
\newcommand{\R}{{\bf{R}}}
\renewcommand{\L}{{\bf{L}}}

\newcommand{\Stab}{\operatorname{Stab}}
\newcommand{\td}{\operatorname{Td}}

\newcommand{\Amp}{\operatorname{Amp}}

\renewcommand{\P}{\mathcal{P}}
\newcommand{\ZZ}{\mathcal{Z}}

\newcommand{\SL}{\operatorname{SL}}

\newcommand{\N}{\operatorname{\mathcal{N}}}
\newcommand{\QQ}{\mathbb{Q}}

\newcommand{\I}{\mathscr{I}}

\begin{document}

\begin{abstract}
We discuss derived categories of coherent sheaves
on algebraic varieties. We focus on the case of non-singular Calabi-Yau varieties and consider two unsolved problems: proving that birational varieties have equivalent derived categories,
and computing the group of derived autoequivalences. We also introduce the space of stability conditions on a triangulated category and explain its relevance to these two problems.
\end{abstract}

\maketitle

\section{Introduction}

In the usual approach to the study of algebraic varieties one focuses directly on
geometric properties of the varieties in question. Thus  one considers embedded curves, hyperplane sections, branched covers
 and so on. A more algebraic approach is to study the varieties indirectly via their (derived) categories of coherent sheaves. 
This second approach has been taken up by an increasing number of  researchers in the last few years. We can perhaps identify three reasons for this new emphasis on categorical methods.

Firstly, algebraic geometers have been attempting to understand string theory. The conformal field theory associated to a variety in string theory contains a huge amount of non-trivial information. However this information seems to be packaged in a categorical way rather than in directly geometric terms.  This perhaps first became clear in Kontsevich's famous homological mirror symmetry conjecture \cite{Kon}. Building on Kontsevich's work, it is now understood that the derived category of coherent sheaves appears in string theory as the category of branes in a topological twist of the sigma model. This means that many dualities in string theory can be described mathematically as equivalences of derived categories.

A second motivation to study varieties via their sheaves is that this approach is expected to generalise more easily to non-commutative varieties.
Although the definition of such objects is not yet clear, there are many interesting examples. In general non-commutative objects have no points in the usual sense, so that direct geometrical methods do not apply. Nonetheless one can hope that the (derived) category of coherent sheaves is well-defined and has  similar properties to the corresponding object in the commutative case.

A third reason is that recent results leave the impression that categorical methods enable one to obtain  a truer description of certain varieties than current geometric techniques allow. For example many equivalences relating the derived categories of pairs of varieties are now known to exist. Any such equivalence points to a close relationship between the two varieties in question, and these relationships are often impossible to describe by other methods. Similarly, some varieties have been found to have interesting groups of derived autoequivalences, implying the existence of symmetries associated to the variety that are not visible in the geometry. It will be interesting to see whether a more geometric framework for these phenomena emerges over the next few years.

Despite a lot of recent work
our understanding of derived categories of coherent sheaves is still quite primitive.
To illustrate this we shall focus in this paper on two easily stated problems, both of which have been around since Bondal and Orlov's pioneering work in the nineties \cite{BO1}, and both of which remain largely unsolved.
Of course the choice of these problems is very much a matter of personal taste. 
Other  survey articles are available by Bondal and Orlov \cite{BO3},  by Hille and van den Bergh \cite{HV} and by Rouquier \cite{Ro}.

Of the many important contributions not discussed in detail  here, we must at least mention Kuznetsov's beautiful work on semi-orthogonal decompositions for Fano varieties.
As well as specific results on Fano threefolds, Kuznetsov has obtained a general theory of homological projective duality with many  applications to derived categories of varieties of interest in classical algebraic geometry. The theory also leads to some extremely interesting derived equivalences between non-commutative varieties.
For more details we refer the reader to the original papers \cite{Ku1,Ku2,Ku3,Ku4}.

\subsection*{Acknowledgements}
I'm very grateful to Antony Maciocia for teaching me about derived categories in the first place, and to Alastair King, Miles Reid and Richard Thomas
for explaining many things since then.
The material in Section 5 has benefitted from discussions with Yukinobu Toda. The author is supported by a Royal Society University Research Fellowship.

\subsection*{Notation}
We write $\D(X)=\operatorname{D^b Coh}(X)$ for the bounded derived category of coherent sheaves on a
variety $X$. All varieties will be assumed to be over the complex numbers. A smooth variety $X$ will be called Calabi-Yau if its canonical bundle $\omega_X$ is trivial. A birational map $f\colon X\to Y$ is called crepant if it is relatively Calabi-Yau, meaning that $\OO_X$ is a relative dualizing object. If $X$ and $Y$ are Gorenstein this is equivalent to the condition $\omega_x=f^*(\omega_Y)$.

\section{Some basic problems}

In this section we describe the basic problems we shall focus on and review some of the known results relating to them.

\subsection{}
The following conjecture was first made by Bondal and Orlov \cite{BO1}.

\begin{conj}[Bondal, Orlov]
\label{conj1}
If $X_1$ and $X_2$ are birational smooth projective Calabi-Yau  varieties of dimension $n$ then there is an equivalence
of categories
$\D(X_1)\to \D(X_2)$.
\end{conj}

So far this is only known to hold in dimension $n=3$. We shall describe the proof in this case in Section 3 below.
There are also some local results available in all dimensions.

\begin{thm}[Bondal, Orlov \cite{BO1}]
\label{BO}
Suppose that $X_1$ and $X_2$ are related by a standard flop, so that there is a diagram
\[
\xymatrix@C=.5em{ &&Y\ar[dll]_{f_1} \ar[drr]^{f_2}\\
X_1&&&& X_2
}
\]
with $f_i$ being the blow-up of a smooth subvariety $E_i\isom\PP^n\subset X_i$ with normal bundle $\OO(-1)^{n+1}$. Then the functor
\[\R f_{2,*}\circ \L f_1^*\colon \D(X_1)\to\D(X_2)\] is an equivalence.
\end{thm}

\begin{thm}[Kawamata \cite{Ka2}, Namikawa \cite{Na1}]
\label{KN}
Suppose that $X_1$ and $X_2$ are related by a Mukai flop, so that there is a diagram
\[
\xymatrix@C=.5em{ &&Y\ar[dll]_{f_1} \ar[drr]^{f_2}\\
X_1\ar[drr]_{g_1}&&&& X_2\ar[dll]^{g_2} \\
&& Z}
\]
with $f_i$ being the blow-up of a smooth subvariety $E_i\isom\PP^n\subset X_i$ with normal bundle equal to the cotangent bundle $\Omega^1(\PP^n)$, and $g_i$ the contraction of $E_i$ to a point. Let $p_i\colon X_1\times_Z X_2\to X_i$ be the projection of the fibre product onto its factors. Then the functor
\[\R p_{2,*}\circ \L p_1^*\colon \D(X_1)\to\D(X_2)\] is an equivalence, whilst in general the functor $\R f_{2,*}\circ \L f_1^*$ is not.
\end{thm}

A result of Wierzba and Wisniewski \cite{WW} states that birationally equivalent complex symplectic fourfolds are related by a finite chain of Mukai flops, so
Theorem \ref{KN} gives

\begin{cor}
Conjecture \ref{conj1} holds when $X_1$ and $X_2$ are complex symplectic fourfolds.
\end{cor}

Kawamata and Namikawa have also studied certain stratified Mukai flops \cite{Ka5,Na2}. In these cases neither a common resolution nor the fibre product give rise to an equivalence, although equivalences do nonetheless exist.
Thus we see that the first problem faced in attempting to prove Conjecture \ref{conj1} is writing down a candidate equivalence.

\subsection{}
A famous variant of Conjecture \ref{conj1} is the derived McKay correspondence, first stated by Reid \cite{Re1,Re2}.

\begin{conj}
\label{conj2}
Let $G\subset \SL(n,\C)$ be a finite subgroup. Let $\D_G(\C^n)$ denote the bounded derived category of the category of $G$-equivariant coherent sheaves on $\C^n$.
Suppose $f\colon Y\to\C^n/G$
is a crepant resolution of singularities.
Then  there is an equivalence of categories $\D(Y)\to \D_G(\C^n)$. 
\end{conj}

Note that $f$ being crepant is equivalent to $Y$ being Calabi-Yau. The condition $G\subset\SL(n,\C)$ is the statement that the canonical bundle of $\C^n$ is $G$-equivariantly trivial. It seems certain that any method which proves Conjecture \ref{conj1} will also apply to this case. On the other hand the McKay correspondence lends itself to more algebraic methods of attack. 

Conjecture \ref{conj2} is known to hold in dimension $n=3$ due to work of Bridgeland, King and Reid \cite{BKR} together with the results of \cite{Br1}. In this case
the quotient variety $\C^3/G$ always has a crepant resolution, a distinguished choice being given by Nakamura's G-Hilbert scheme.
This scheme is perhaps best thought of as a moduli space of representations of the skew group algebra $\C[x,y,z]*G$
that are stable with respect to a certain choice of stability condition. Craw and Ishii \cite{CI} showed how to obtain  other resolutions  by varying the stability condition. 

The conjecture  has recently been shown to hold in two other situations, although the methods of proof are very different from that employed in the dimension three case. 
Firstly, Bezrukavnikov and Kaledin \cite{BK} used characteristic $p$ methods and deformation quantization to deal with symplectic group quotients.

\begin{thm}[Bezrukavnikov, Kaledin \cite{BK}] Conjecture \ref{conj2} holds when $G$ preserves
a complex symplectic form on $\C^n$.
\end{thm}

Secondly, a special case of Kawamata's work \cite{Ka3} on the effect or toroidal flips and flops on derived categories is the following result. 

\begin{thm}[Kawamata \cite{Ka4}]
\label{abmc}
Conjecture \ref{conj2} holds when $G$ is abelian.
\end{thm}

Note that in this case all resolutions are toric. Kawamata's proof uses the toric minimal model programme together with explicit computations of $\Hom$ spaces.

\subsection{}
The second problem we wish to focus on is even easier to state, namely

\begin{prob}
\label{prob}
Compute the group of autoequivalences of the derived category of a smooth projective Calabi-Yau variety $X$.
\end{prob}

If $X$ is any smooth projective variety there are certain standard autoequivalences of $\D(X)$, namely the automorphisms of $X$ itself, twists by elements of $\Pic(X)$, and the shift functor $[1]$. Bondal and Orlov \cite{BO1,BO2} showed that if $X$ has ample canonical or anticanonical bundle, then these standard equivalences generate  $\Aut\D(X)$. Thus the group $\Aut\D(X)$ is not so interesting in this case. In fact the group $\Aut\D(X)$ seems to be most interesting when $X$ is Calabi-Yau, so we concentrate on this case.

Certain types of objects are known to define autoequivalences. The most well-known are  spherical objects.

\begin{thm}[Seidel, Thomas \cite{ST}]
Let $X$ be a smooth projective Calabi-Yau of dimension $n$ and let $S\in\D(X)$ be a spherical object, which is to say an object satisfying
\[
\Hom^k_{\D(X)}(S,S)=\biggl\{\begin{array}{ll} \C & \text{ if }k=0 \text{ or }n,
\\ 0 &\text{ otherwise.}\end{array}\biggr.
\]
Then there is an autoequivalence $\Phi_S\in\Aut(\D)$ such that for any object $E\in\D$ there is an exact triangle
\[\Hom_{\D(X)}(S,E)\tensor S\lra E\lra \Phi_S(E).\]
\end{thm}

Szendr{\H o}i \cite{Sz} showed how certain configurations of surfaces in Calabi-Yau threefolds lead to interesting subgroups of autoequivalences of the corresponding derived categories.
More recently, Huybrechts and Thomas \cite{HT} have studied projective space objects, which also give rise to autoequivalences.

Problem \ref{prob} is already  non-trivial when $X$ is an elliptic curve.
Since $X$ can be identified with its dual abelian variety $\Pic^0(X)$ in the dimension one case, the original Fourier transform
of Mukai \cite{Mu1} defines an autoequivalence $\F\in\Aut\D(X)$.
This enables one to prove

\begin{thm}
\label{ell}
The group $\Aut\D(X)$ is generated by the standard autoequivalences together with $\F$. There is an exact sequence
\[\begin{CD}
1@>>> \Z\times (\Aut(X)\ltimes\Pic^0(X))@>>> \Aut\D(X)@>f>> \SL(2,\Z)\lra 1.
\end{CD}
\]
where the factor of $\Z$ is generated by the double shift $[2]$, and the map $f$ is defined by
\[\mat{\rk(\Phi(E))}{\dg(\Phi(E))}=f(\Phi)\mat{\rk(E)}{\dg(E)}\]
for all $E\in\D(X)$.
\end{thm}

Choosing a degree one line bundle $L$ on $X$ gives a splitting of the exact sequence in Theorem \ref{ell}, because the equivalences $-\tensor L$ and $\F$ generate a subgroup of $\Aut\D(X)$ isomorphic to $\SL(2,\Z)$. This $\SL(2,\Z)$ action on $\D(X)$ discovered by Mukai was perhaps an early pointer to homological mirror symmetry.

Therorem \ref{ell}  was generalised by Orlov \cite{Or} who calculated the group $\Aut\D(X)$ for all abelian varieties $X$. 
The group of derived autoequivalences of the minimal resolution of a Kleinian singularity was studied in \cite{Br4,IU}. But beyond these results not much is known. We shall consider the case when $X$ is a K3 surface in Section 6.

\section{Threefold flops}

In this section we shall be concerned with the following result.

\begin{thm}
\label{mine}
Conjecture \ref{conj1} holds when $X_1$ and $X_2$ have dimension three.
\end{thm}

Of course in dimensions two or less, birational Calabi-Yau varieties are isomorphic.
We start by giving a bare outline of the original proof of Theorem \ref{mine}
and then return to make some further remarks on each step. Full details of the proof can be found in \cite{Br1}.
A more direct proof has since been found by Kawamata \cite{Ka2} using results of Van den Bergh \cite{VdB}. Nonetheless, the original proof is worth describing since it gives a moduli-theoretic interpretation of threefold flops which is interesting in its own right.

\medskip
\noindent{\bf Sketch of proof}.
The first step is to apply the minimal model programme in the shape of a result of Koll{\'a}r \cite{Kol} which says that
any birational transform between smooth projective threefolds that preserves the canonical class can be split into a finite sequence of flops.
Thus we can reduce to the case where we have a flopping diagram
\[
\xymatrix@C=.5em{ Y\ar[drr]_f&&&& W\ar[dll]^g \\
&&X
}
\]
Here $Y$ and $W$ are smooth and projective, $X$ has Gorenstein terminal singularities, and $f$ and $g$ are crepant birational maps contracting only
finitely many rational curves.

The second step is to define a full subcategory $\Per\subset \D(Y)$ consisting of objects $E\in\D(Y)$ satisfying
the following three conditions

\begin{itemize}
\item[(a)] $H^i(E)=0$ unless $i=0$ or $-1$,

\item[(b)] $\R^1 f_*\, H^0(E)=0$ and $\R ^0 f_*\, H^{-1}(E)=0$,

\item[(c)] $\Hom_Y(H^0(E),C)=0$ for any sheaf $C$ on $Y$ satisfying $\R f_*(C)=0$.
\end{itemize}
In fact $\Per$ is the heart of a bounded t-structure on $\D(Y)$ and hence is abelian. 
We call the objects of $\Per$ perverse coherent sheaves.

Define a perverse point sheaf on $Y$ to be an object of $\Per$
which has the same Chern classes as the structure sheaf of a point $y\in Y$, and which is a quotient of $\OO_Y$ in $\Per$.
The third step is to show that there is a fine moduli space $M$ for perverse point sheaves on $Y$. 
Since perverse point sheaves are not sheaves in general, but complexes of sheaves with nontrivial cohomology
in more than one place, the usual moduli space techniques do not apply. In fact we can side-step this issue as follows.

By definition each perverse point sheaf $E$ fits into a short exact sequence
\[0\lra F\lra \OO_{Y}\lra E\lra 0\]
 in $\Per$. Taking the long exact sequence in cohomology shows that $F$ is actually a sheaf, although not necessarily torsion-free. One can construct a fine moduli space $M$ for the objects $F$ using a standard GIT argument, and since the $E$s and $F$s determine one another this space is also a fine moduli space for perverse point sheaves.

The moduli space $M$ has a natural map to $X$ since it is easily shown that if $E$ is a perverse point sheaf on $Y$ then $\R f_*(E)$ is the structure sheaf of a point of $X$.
The fourth step is to show that $M$ is smooth  and that the universal family $\P$ on $M\times Y$ induces a Fourier-Mukai equivalence $\Phi\colon \D(M)\to\D(Y)$. This involves an application of the famous new intersection  theorem (see for example \cite{NIT}).
Once one has this it is easy to see that $M$ with its natural map to $X$ is the flop of $Y\to X$ and hence can be identified with $W$.

Finally, we should note that  Chen \cite[Prop. 4.2]{Ch}  was able to show that the universal family $\P$  is isomorphic to 
the object $\OO_{W\times_X Y}$. Thus the resulting functor $\Phi$ is of the same form used in Bondal and Orlov's result Theorem \ref{BO}. Later Kawamata \cite{Ka2} was able to show directly that this functor gives an equivalence.
\qed

\begin{example}
To understand why the moduli of perverse point sheaves gives the flop consider the simplest example
when $f\colon Y\to X$ is the contraction of a
non-singular rational curve $C$ with normal bundle
$\OO_C(-1)\oplus\OO_C(-1)$.

Structure sheaves of points $y\in Y$
are objects of the
category $\Per$ for all $y\in Y$. 
However, if $y\in C$ then the nonzero map $\OO_C(-1)\to \I_y$ means that $\I_y$ is not an object of $\Per$, so that $\OO_y$ is not a quotient of $\OO_Y$ in $\Per$ and hence is not a perverse point sheaf. 
In fact for $y\in C$, the sheaf
$\OO_y$ fits into the exact sequence
\begin{equation}
\label{a}
0\lra \OO_C(-1)\lra \OO_C\lra \OO_y\lra 0.
\end{equation}
Now $\OO_C$ is a perverse sheaf, but $\OO_C(-1)$ is
not, so that the triangle in $\D(Y)$ arising from (\ref{a})
does not define an exact sequence in
$\Per$. However the complex obtained by shifting $\OO_C(-1)$
to the left by one place \emph{is} a perverse sheaf, so there
is an exact sequence of perverse sheaves
\begin{equation}
\label{b}
0\lra \OO_C\lra\OO_y\lra\OO_C(-1)[1]\lra 0
\end{equation}
which could be thought of as destabilizing $\OO_y$.

Flipping the extension of perverse sheaves (\ref{b}) gives
 objects of $\Per$ fitting into an exact sequence of
perverse sheaves
\begin{equation}
\label{c}
0\lra \OO_C(-1)[1]\lra E\lra \OO_C\lra 0.
\end{equation}
It is easy to see that these objects $E$ are perverse point sheaves. Note that they are
not sheaves, indeed any such object has two nonzero
cohomology sheaves $H^{-1}(E)=\OO_C(-1)$ and $H^0(E)=\OO_C$.
 Roughly speaking, the space $W$ is obtained from $X$
by replacing the rational curve $C$ parameterising extensions
(\ref{b}) by another rational curve $C'$ parameterising
extensions (\ref{c}).
\end{example}
\medskip

We now make some further remarks on the proof with an eye to generalising the method to higher dimensions, although, as we shall see, the prospects for doing this do not seem particularly good.

It is clear that the key step in the proof is the introduction  of the category $\Per$. There are two possible ways to arrive at this definition as we shall now explain.
Unfortunately neither of them seems to generalise directly to higher dimensional situations.

\smallskip
\noindent{\bf Construction 1}. The first approach is to use the theory of tilting in abelian categories as developed by Happel, Reiten and Smal\o \cite{HRS}. This is an abstraction of the notion of tilting of algebras due to Brenner and Butler \cite{BB}. In fact in  our situation the full subcategory
\[\F=\{E\in\Coh(Y):f_*(E)=0\}\]
is the torsion-free subcategory of a torsion pair on $\Coh(Y)$.
The corresponding torsion subcategory consists of objects $E\in\Coh(Y)$ for which $\R^1 f_*(E)=0$ and $\Hom_Y(E,C)=0$ for any sheaf $C$ on $Y$ satisfying $\R f_*(C)=0$. Tilting with respect to this torsion pair gives the t-structure on $\D(Y)$ whose heart is $\Per$.

\smallskip
\noindent{\bf Construction 2}. Van den Bergh \cite{VdB} found another characterisation  of $\Per$ involving sheaves of non-commutative algebras. He first introduced a particular locally-free sheaf $P$ on $Y$ and considered the coherent sheaf of algebras $\A=\R f_* \lEnd(P)$ on $X$. Then he showed that there is an equivalence
\[\R f_*\lHom_{\OO_Y}(P,-)\colon \D(Y)\lra \D(\A),\]
where $\D(\A)$ is the bounded derived category of the abelian category $\Coh(\A)$  of coherent sheaves of $\A$-modules on $X$. Pulling back the standard t-structure on $\D(\A)$ gives the t-structure on $\D(Y)$ whose heart is $\Per$. In particular the above derived equivalence restricts to give an equivalence of abelian subcategories $\Per\to \Coh(\A)$. 

\smallskip

The third step in  the proof of Theorem \ref{mine} was to construct a fine moduli space for perverse point sheaves. The trick we used in dimension three will not work in higher dimensions. Instead one must construct the moduli space directly. Recent results of Inaba \cite{In} and Toen and Vaquie \cite{To} construct spaces parameterising very general classes of objects of derived categories. For applications we need a more geometric approach. We first need to understand what it means for an object of a triangulated category to be stable. The general notion of a stability condition
introduced in the next section sheds some light on this. After that we need to construct good
 moduli spaces of stable objects, where good might mean for example that the moduli space is a projective scheme.  
We refer to Abramovich and Polishchuk \cite{AP} for some recent progress in this direction.

The final step in  the proof of Theorem \ref{mine} was very much dependent on the dimension three assumption. The estimates needed to apply the intersection theorem just do not work in  higher dimensions. In fact, in dimension four and higher, one can obtain singular varieties by flopping non-singular ones, so we cannot expect to prove that the relevant moduli space is smooth. Derived categories of singular varieties are still poorly understood at present, a point to which we shall return in the last section.

\section{Stability conditions}

The notion of a stability condition was introduced in \cite{Br2} as a way to understand Douglas' work on $\pi$-stability for D-branes in string theory \cite{Do}. Here we wish to emphasise the purely mathematical aspects of this definition. For more on the connections with string theory see \cite{Br5}.

In the context of the present article stability conditions are relevant for three reasons. Firstly, the choice of a stability condition picks out classes of stable objects for which one can hope to form well-behaved moduli spaces. Secondly the space of all stability conditions $\Stab(\D)$ allows one to bring geometric methods to bear on the problem of understanding t-structures on $\D$. Finally, the space $\Stab(\D)$ provides a complex manifold on which the group $\Aut(\D)$  naturally acts. 

Throughout this section $\D$ is a fixed triangulated category. We shall assume that $\D$ is essentially small, i.e. equivalent to a category whose objects form a set.
The Grothendieck group of $\D$ is denoted $K(\D)$. For details on the results of this section see \cite{Br2}.

\subsection{}
The definition of a stability condition is as follows.

\begin{defn}
\label{pemb}
A stability condition $\sigma=(Z,\P)$ on $\D$
consists of
a group homomorphism
$Z\colon K(\D)\to\C$ called the central charge,
and full additive
subcategories $\P(\phi)\subset\D$ for each $\phi\in\R$,
satisfying the following axioms:
\begin{itemize}
\item[(a)] if $E\in \P(\phi)$ then $Z(E)=m(E)\exp(i\pi\phi)$ for some
 $m(E)\in\R_{>0}$,
\item[(b)] for all $\phi\in\R$, $\P(\phi+1)=\P(\phi)[1]$,
\item[(c)] if $\phi_1>\phi_2$ and $A_j\in\P(\phi_j)$ then $\Hom_{\D}(A_1,A_2)=0$,
\item[(d)] for each nonzero object $E\in\D$ there is a finite sequence of real
numbers
\[\phi_1>\phi_2> \cdots >\phi_n\]
and a collection of triangles
\[
\xymatrix@C=.5em{
0_{\ } \ar@{=}[r] & E_0 \ar[rrrr] &&&& E_1 \ar[rrrr] \ar[dll] &&&& E_2
\ar[rr] \ar[dll] && \ldots \ar[rr] && E_{n-1}
\ar[rrrr] &&&& E_n \ar[dll] \ar@{=}[r] &  E_{\ } \\
&&& A_1 \ar@{-->}[ull] &&&& A_2 \ar@{-->}[ull] &&&&&&&& A_n \ar@{-->}[ull] 
}
\]
with $A_j\in\P(\phi_j)$ for all $j$.
\end{itemize}
\end{defn}

Given a stability condition $\sigma=(Z,\P)$ as in the definition, each subcategory
$\P(\phi)$ is abelian. The nonzero objects of $\P(\phi)$ are said to be
{semistable of phase $\phi$ in $\sigma$}, and the simple objects
of $\P(\phi)$ are said to be stable.
It follows from the other axioms that the decomposition of
an
object $0\neq E\in\D$ given by
axiom (d) is uniquely defined up to isomorphism.
Write
$\phi^+_{\sigma}(E)=\phi_1$ and
$\phi^-_{\sigma}(E)=\phi_n$. 
The {mass} of $E$ is defined to be
the positive real number
$m_{\sigma}(E)=\sum_i |Z(A_i)|$.

For any interval $I\subset\RR$, define $\P(I)$ to be the extension-closed
subcategory of $\D$ generated by the subcategories $\P(\phi)$ for
$\phi\in I$. Thus, for example, the full subcategory $\P((a,b))$
consists of the zero objects of $\D$ together with those
objects $0\neq E\in\D$ which satisfy
$a<\phi_\sigma^-(E)\leq\phi_\sigma^+(E)<b$.
A stability condition is called {locally finite} if
there is some $\epsilon>0$ such that each quasi-abelian category
$\P((\phi-\epsilon,\phi+\epsilon))$ is of finite length.

\subsection{}
If $\sigma=(Z,\P)$ is a stability condition on $\D$, one can check that the full subcategory $\A=\P((0,1])$ is the heart of a bounded t-structure on $\D$.
We call it the heart of $\sigma$. The stability condition $\sigma$ can easily be reconstructed from its heart $\A$ together with the central charge map $Z$. 
Since this will be important in Section 5 we spell it out in a little more detail.

Define a   stability function on
an abelian  category $\A$ to be a group homomorphism $Z\colon K(\A)\to\C$
such that
\[0\neq E\in\A \implies Z(E)\in\R_{>0}\,\exp({i\pi\phi(E)})\text{ with
}0<\phi(E)\leq 1.\]
The real number $\phi(E)\in(0,1]$ is called the phase of the object $E$.

A nonzero object $E\in\A$ is said to be
semistable
with respect to $Z$
if every subobject $0\neq A\subset E$ satisfies $\phi(A)\leq\phi(E)$.
The  stability function $Z$ is said to have the Harder-Narasimhan property if
every nonzero object $E\in\A$ has
a finite filtration
\[0=E_0\subset E_1\subset \cdots\subset E_{n-1}\subset E_n=E\]
whose factors $F_j=E_j/E_{j-1}$ are semistable objects of $\A$ with
\[\phi(F_1)>\phi(F_2)>\cdots>\phi(F_n).\]

Given a stability condition on $\D$, the central charge defines a  stability function on
its heart  $\A=\P((0,1])\subset\D$, and the decompositions of axiom (d) give Harder-Narasimhan filtrations.
Conversely, given a bounded t-structure on $\D$ together with a stability function $Z$ on its heart $\A\subset\D$, we can define subcategories
$\P(\phi)\subset\A\subset\D$ to be the semistable objects in $\A$ of phase $\phi$ for each $\phi\in(0,1]$. Axiom (b) then fixes $\P(\phi)$ for all $\phi\in\R$.
Thus we have

\begin{prop}
\label{pg}
To give a stability condition on $\D$ is
equivalent to giving a bounded t-structure on $\D$ and a 
 stability function on its heart with the Harder-Narasimhan property.
\end{prop}

\subsection{}
The set $\Stab(\D)$ of locally-finite stability conditions on $\D$ has a natural topology induced by the
metric
\[d(\sigma_1,\sigma_2)=\sup_{0\neq E\in\D}
\bigg\{|\phi^-_{\sigma_2}(E)-\phi^-_{\sigma_1}(E)|,|\phi^+_{\sigma_2}(E)
-\phi^+_{\sigma_1}(E)|
,|\log \frac{m_{\sigma_2}(E)}{m_{\sigma_1}(E)}|\bigg\}
\in[0,\infty].\]
The following  result was proved in \cite{Br1}. Its slogan is that deformations of the central charge
lift uniquely to deformations of the stability condition.

\begin{thm}
\label{lasty}
For each connected component $\Sigma\subset\Stab(\D)$
there is a linear subspace $V(\Sigma)\subset \Hom_{\Z}(K(\D),\Z)$
with a well-defined linear topology and a local homeomorphism
$\ZZ\colon\Sigma\to V(\Sigma)$ which maps a stability condition $(Z,\P)$
to its central charge $Z$.
\end{thm}

If $X$ is a smooth projective variety we write $\Stab(X)$ for the set of locally-finite stability conditions on $\D(X)$ for which
the central charge $Z$ factors via the Chern character map
$\ch\colon K(X)\to H^*(X,\QQ)$.
Theorem \ref{lasty} immediately implies that $\Stab(X)$ is a finite-dimensional complex manifold.
By Proposition \ref{pg}, the points of $\Stab(X)$ parameterise bounded t-structures on $\D(X)$, together with the extra data of the map $Z$.

So far the manifolds $\Stab(X)$ have only been computed for varieties $X$ of dimension one \cite{Ma,Ok}. The case of K3 and abelian surfaces was studied in \cite{Br3}; we shall return to the K3 case in Section 6. Various non-compact examples have also been considered (see e.g. \cite{Br4,Toda}). 

\section{Stability conditions and threefold flops}

Let $f\colon Y\to X$ be a small, crepant birational map of threefolds with $Y$ smooth. Thus $X$ has terminal Gorenstein singularities
and $f$ contracts a finite number of rational curves $C$, each satisfying $K_Y\cdot C=0$.  The aim of this section is to show how the category of perverse sheaves defined in Section 3 arises naturally by considering stability conditions on 
$Y$. We shall go into more detail than we have in previous sections, since the results provide an excellent example of how stability conditions can be used to relate the derived category to geometry. For full proofs the reader can consult Toda's paper \cite{Toda}. The string theory point of view on the same material is described in \cite{As}. 

\subsection{}
Define $\D(Y/X)$ to be the full subcategory of $\D(Y)$
consisting of objects supported on (a fat neighbourhood of) the exceptional locus of $f$. 
Since $\D(Y/X)$ only depends on a formal neighbourhood of the singular locus of $X$ we may as well assume that $X$ is the spectrum of a complete local ring.

One can easily show that there is an isomorphism of abelian groups
\[K(\D(Y/X))\isom N_1(Y/X)\oplus\Z.\]
It will be convenient to consider the codimension one slice of the full space $\Stab(\D(Y/X))$ consisting of stability conditions satisfying the
additional condition that $Z(\OO_y)=-1$, where $\OO_y$ is the structure sheaf of a point $y\in Y$. We shall denote it simply by $\Stab(Y/X)$.
Since the dual of $N_1(Y/X)$ is naturally $N^1(Y/X)$, the central charge $Z\colon K(\D(Y/X))\to\C$ of a stability condition in $\Stab(Y/X)$ is defined by an element $\beta+i\omega$ of the vector space $N^1(Y/X)\tensor\C$. Explicitly the correspondence is given by
\[\tag{$\dagger$} Z(E)=(\beta+i\omega)\cdot \ch_2(E) -\ch_3(E).\]
The map $\ZZ$ of Theorem \ref{lasty} now becomes a map
\[\ZZ\colon\Stab(Y/X)\lra N^1(X/Y)\tensor\C.\]

The standard t-structure on $\D(Y)$ induces a bounded t-structure on $\D(Y/X)$ whose heart $\Coh(Y/X)=\Coh(Y)\cap\D(Y/X)$ is the subcategory of $\Coh(Y)$ consisting of sheaves supported on the exceptional locus. A simple application of Proposition \ref{pg} gives

\begin{prop}
\label{bored}
A  stability function for the t-structure with heart $\Coh(Y/X)\subset\D(Y/X)$ is given by $(\dagger)$ with  $\beta,\omega\in N^1(Y/X)\tensor{\R}$ such that $\omega\in\Amp(Y/X)$ lies in the ample cone. This stability function has the Harder-Narasimhan property.
Thus there is a region in $U(Y/X)\subset\Stab(Y/X)$ isomorphic to the complexified ample cone of $Y/X$.
\end{prop}

\subsection{}
Suppose now that $\sigma=(Z,\P)$ is a stability condition in the boundary of the region $U(Y/X)$.
Then $Z$ is given by the formula above for some $\beta,\omega\in N^1(Y/X)\tensor{\R}$ with $\omega$ an $f$-nef divisor which is not ample.
There is a unique birational map $f'\colon Y\to Y'$ factoring $f\colon Y\to X$ such that $\omega=(f')^*(\omega')$ with $\omega'$ an ample $\R$-divisor on $Y'$.
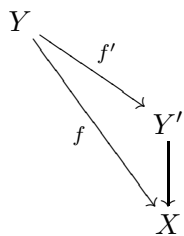
\[
\xymatrix@C=1.5em{ Y\ar[drr]^{f'}\ar[ddrr]_f \\
&&Y' \ar[d]\\
&&X
}
\]
The map $f'$ is also crepant and contracts precisely those irreducible rational curves $C\subset Y$ satisfying $\omega\cdot C=0$.

Let $C$ be a rational curve contracted by $f'$. Note that for all stability conditions in $U(Y/X)$, the objects $\OO_C(k)$  are stable, since their only proper quotients in $\Coh(Y/X)$ are skyscraper sheaves. It follows that these objects are at least semistable in $\sigma$, and hence $Z(\OO_C(k))\neq 0$. Since $\omega\cdot C=0$ we must have $\beta\cdot C\notin\Z$.

The map $Z$ defined by $(\dagger)$  no longer defines a  stability function for $\Coh(Y/X)$ because if $C$ is contracted by $f'$ and $k\ll 0$ then $Z(\OO_C(k))$ lies on the positive real axis. Thus the heart of the stability condition $\sigma$ cannot be $\Coh(Y/X)\subset \D(Y/X)$.

Instead, consider the perverse t-structure on $\D(Y)$ induced by the contraction $f'$ and restrict it to $\D(Y/X)$ to give a bounded t-structure whose heart $\Perw\cap \D(Y/X)$ we shall also denote $\Perw$. One can then show that providing $0<\beta\cdot C<1$ for all irreducible curves $C$ with $\omega\cdot C=0$ then the heart of $\sigma$ is indeed equal to $\Perw$.

Suppose for a moment that $\sigma$ lies on a codimension one face of the boundary of $U(Y/X)$, by which we mean that $\Pic(Y/Y')=\Z$. Note that $\Pic(Y/X)$ acts on $\Stab(Y/X)$ and on the region $U(Y/X)$. Thus twisting by a line bundle we can always assume that the condition $0<\beta\cdot C<1$ holds.

Conversely, considering the definition of $\Perw$ as a tilt, one can show that if the above condition on $\beta$ and $\omega$ holds then the map $Z$ defined by $(\dagger)$  gives a  stability function on $\Perw$, and this is easily checked to satisfy the Harder-Narasimhan property.
These arguments give

\begin{prop}
\label{fish}
Let $f'\colon Y\to Y'$ be a birational map factoring $f\colon Y\to X$.
A  stability function for the t-structure with heart $\Perw$ is given by $(\dagger)$
for classes $\beta,\omega\in N^1(Y/X)\tensor{\R}$ such that
$\omega=(f')^*(\omega')$ for some $\omega'\in\Amp(Y'/X)$ and $0<\beta\cdot C<1$
for all curves $C$ contracted by $f'$. The resulting stability conditions lie in the boundary of the region $U(Y/X)$.
Conversely, any stability condition lying on a codimension one face of the boundary of $U(Y/X)$ is of this form up to a twist by a line bundle.
\end{prop}
The key point we wish to emphasize is that the perverse t-structure arises naturally by passing to the boundary of the region  $U(Y/X)$.

\subsection{}
Consider now an open neighbourhood of a point $\sigma$ lying in the boundary of $U(Y/X)$. This region will contain stability conditions with central charges $Z$ defined by $(\dagger)$  with $\omega$ lying in the ample cone of a different birational model of $f\colon Y\to X$. It can be shown that all these stability conditions can be obtained by pulling back stability conditions of the form given in Proposition \ref{bored} by derived equivalences.

For example, suppose $g\colon  W\to X$ is obtained from $f$ by taking the flop of a contraction $f'\colon  Y\to Y'$ factoring $f$ with $\Pic(Y/Y')\isom\Z^{\oplus d}$. 
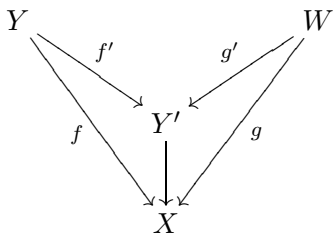
\[
\xymatrix@C=1.5em{ Y\ar[drr]^{f'}\ar[ddrr]_f&&&& W\ar[dll]_{g'}\ar[ddll]^g \\
&&Y' \ar[d]\\
&&X
}
\]
Such a map $f'$ defines a codimension $d$ face of the ample cone $\Amp(Y/X)$.
Let $\Phi\colon \D(Y)\lra \D(W)$ be the equivalence with kernel $\OO_{Y\times_{Y'} W}$ on $Y\times W$.
Then $\Phi$ induces an isomorphism fitting into a diagram
\[
\xymatrix@C=.8em{ \Stab(Y/X)\ar[d]_{\ZZ}\ar[rrrr]^{\Phi} &&&& \Stab(W/X)\ar[d]_{\ZZ}\\
 N^1(Y/X)\tensor\C\ar[rrrr]^{(g^{-1}\circ f)_*} &&&& N^1(W/X) \tensor\C
}
\]
where the isomorphism $N^1(Y/X)\to N^1(W/X)$ is induced by the codimension one isomorphism $g^{-1}\circ f\colon Y\dashrightarrow W$. One can show that the closure of the inverse image $\Phi^{-1}(U(W/X))$ 
intersects the closure of $U(Y/X)$ along a real codimension $d$  component of the boundary.

It is well known that the general hyperplane section of $f\colon Y\to X$ is the minimal resolution of some Kleinian ADE surface singularity.\footnote{Note added in 2019. In fact this is not generally true and should be added as an assumption in Theorem \ref{tda}. Compare \cite[Section 3]{Toda}. For a description of the  space of stability conditions when this condition fails see  Y. Hirano and M. Wemyss,  `Stability conditions for 3-fold flops', arXiv:1907:09742.} Furthermore the abelian group $N_1(Y/X)$ can be identified with the root lattice of the corresponding simple Lie algebra. Let $\Lambda\subset N_1(Y/X)$ be the set of roots.

\begin{thm}[Toda, \cite{Toda}]\label{tda}
The map \[\ZZ\colon\Stab(Y/X)\lra N^1(Y/X)\tensor {\C}\]
restricted to the connected component of $\Stab(Y/X)$ containing $U(Y/X)$ is a covering map
of the hyperplane complement
\[\{\beta+i\omega\in N^1(Y/X)\tensor\C: \beta+i\omega\cdot C\notin\Z \text{ for all }C\in\Lambda\}.\]
For each smooth birational model $W\to X$ and each Fourier-Mukai equivalence
\[\Phi\colon \D(Y)\lra \D(W)\]
over $X$ there is a region $\Phi^{-1}(U(W/X))\subset \Stab(Y/X)$
isomorphic to the complexified ample cone of $W/X$.  The closures of these regions cover a connected component of $\Stab(Y/X)$, and any two of the regions  are either disjoint or equal.
\end{thm}

\section{Stability conditions on K3 surfaces}

Suppose that $X$ is an algebraic K3 surface over $\C$.
In this section we consider stability conditions on $\D(X)$ and explain
how they may help to determine the group of autoequivalences $\Aut\D(X)$.
Full details can be found in \cite{Br3}.

\subsection{}
Following Mukai \cite{Mu2}, one introduces the
extended cohomology lattice of $X$ by using the formula
\[\big((r_1,D_1,s_1),(r_2,D_2,s_2)\big)=D_1\cdot D_2-r_1 s_2-r_2 s_1\]
to define a symmetric bilinear form on the
cohomology ring
\[H^{*}(X,\Z)=H^0(X,\Z)\oplus H^2(X,\Z)\oplus H^4(X,\Z).\]
The resulting lattice $H^*(X,\Z)$ is even and non-degenerate
and has  signature $(4,20)$.
Let $H^{2,0}(X)\subset H^2(X,\C)$
denote the one-dimensional complex subspace spanned by the
class of a nonzero holomorphic two-form $\Omega$ on $X$.
An isometry \[\varphi\colon H^*(X,\Z)\to H^*(X,\Z)\] is called a Hodge
isometry if $\varphi\tensor\C$ preserves this subspace.
The group of Hodge isometries of
$H^*(X,\Z)$ will be denoted $\Aut H^*(X,\Z)$.

The Mukai vector of an object $E\in\D(X)$ is the element of the
sublattice
\[\N(X)=\Z\oplus\NS(X)\oplus\Z=H^*(X,\Z)\cap
\Omega^{\perp}\subset
H^*(X,\C)\]
defined by the formula
\[v(E)=(r(E),c_1(E), s(E))=\ch(E)\sqrt{\td(X)}\in H^*(X,\Z),\]
where $\ch(E)$ is the Chern character of $E$ and $s(E)=\ch_2(E)+r(E)$.
The Mukai bilinear form makes $\N(X)$ into
an even lattice of signature $(2,\rho)$, where
$1\leq\rho\leq 20$ is the Picard number of $X$.
The Riemann-Roch theorem shows that this form is
the negative of the Euler form, that is,
for any pair of objects $E$ and $F$ of $\D(X)$
\[\eu(E,F)=\sum_i (-1)^i \dim_{\C} \Hom_X^i(E,F)=-(v(E),v(F)).\]

A result of Orlov \cite{Or},
extending work of Mukai \cite{Mu2}, shows
that every exact autoequivalence of $\D(X)$ induces a Hodge isometry of the
lattice $H^*(X,\Z)$. Thus there is a group homomorphism
\[\varpi\colon \Aut\D(X)\lra\Aut H^*(X,\Z).\]
The kernel of this homomorphism will be denoted $\Aut^0\D(X)$.

For any point $\sigma=(Z,\P)\in\Stab(X)$ the central charge $Z(E)$ depends only
on the Chern character of $E$ and hence can be written in the form
\[Z(E)=(\mho,v(E))\]
for some vector $\mho\in \N(X)\tensor\C$. 
Thus the map $\ZZ$ of Theorem \ref{lasty} becomes a map 
\[\ZZ\colon\Stab(X)\to\N(X)\tensor\C\]
sending a stability condition to the corresponding vector
$\mho\in\N(X)\tensor\C$.

Define an open subset
\[\Pol^{\pm}(X)\subset\N(X)\tensor\C\]
consisting of those vectors which span positive definite two-planes
in $\N(X)\tensor\R$.
This space has two connected components that are exchanged
by complex conjugation.
Let $\Pol^+(X)$ denote the component containing vectors of the form
$\exp({i\omega})=(1,i\omega,-\omega^2/2)$ for ample divisor classes
$\omega\in\NS(X)\tensor\R$.
Set \[\Delta(X)=\{\delta\in\N(X):(\delta,\delta)=-2\}\]  and for each
$\delta\in\Delta(X)$ let
\[\delta^{\perp}=\{\mho\in\N(X)\tensor\C:(\mho,\delta)=0\}\subset
\N(X)\tensor\C\]
be the
corresponding complex hyperplane. The following result was proved in \cite{Br3}.
 
\begin{thm}
\label{first}
There is a connected component $\Sigma(X)\subset\Stab(X)$
that is mapped by $\ZZ$ onto the open subset
\[\Pol_0^+(X)=\Pol^+(X)\setminus\bigcup_{\delta\in\Delta(X)}\delta^{\perp}
\subset\N(X)\tensor\C.\]
Moreover, the induced map
$\ZZ\colon \Sigma(X)\to\Pol^+_0(X)$
is a regular covering map and the subgroup of $\Aut^0 \D(X)$
which preserves the connected component $\Sigma(X)$ acts freely on
$\Sigma(X)$ and is the group of deck transformations of $\ZZ$.
\end{thm}

Unfortunately, Theorem \ref{first} is
not enough to determine the group $\Aut \D(X)$. Nonetheless it provides a hope to solve this problem by
studying the geometry of $\Stab(X)$. 
The obvious thing to conjecture is as follows.

\begin{conj}
\label{conj}
The action of $\Aut \D(X)$ on $\Stab(X)$ preserves the connected
component $\Sigma(X)$. Moreover $\Sigma(X)$ is simply-connected.
\end{conj}

This conjecture would imply
that there is a short exact
sequence of groups
\[1\lra \pi_1 \Pol_0^+(X)\lra
\Aut \D(X)\lRa{\varpi} \Aut^+ H^*(X,\Z)\lra
1\]
where \[\Aut^+ H^*(X,\Z)\subset\Aut H^*(X,\Z)\] is the index 2 subgroup
consisting of elements which do not exchange the two components of
$\Pol^{\pm}(X)$.

As a final remark in this section note that Borcherds' work on modular forms \cite{Bo} allows one to write down product expansions for holomorphic
functions on $\Stab(X)$ that are invariant under the group $\Aut\D(X)$. It would be interesting to connect these formulae with
counting invariants for stable objects in $\D(X)$.

\section{Derived categories and the minimal model programme}

In this section we briefly mention some further recent work on birational geometry and derived categories of coherent sheaves and give some references.

\subsection{}
Conjecture \ref{conj1} is the simplest version of a much more general set of conjectures about derived categories and birational geometry.
The basic idea, due to Bondal and Orlov, is that
each of the operations in the minimal model programme should induce fully faithful embeddings of derived categories. Thus for example, if $f\colon Y\to X$ is a birational map of smooth projective varieties, then $\R f_*(\OO_Y)=\OO_X$, and the projection formula implies that the functor
\[\L f^*\colon \D(X)\lra\D(Y)\]
is full and faithful. 
Much more generally Bondal and Orlov conjectured

\begin{conj}[Bondal, Orlov]
\label{c3}
Suppose
\[
\xymatrix@C=.5em{ &&Y\ar[dll]_{f_1} \ar[drr]^{f_2}\\
X_1&&&& X_2
}
\]
are birational maps of smooth projective varieties. Suppose the divisor 
\[K_{X_2/X_1}=f_2^*(K_{X_2})-f_1^*(K_{X_1})\]
 is effective. Then there is a fully faithful functor $F\colon \D(X_1)\into \D(X_2)$.  If $K_{X_2/X_1}=0$ then  $F$ is an equivalence.
\end{conj}

One might well imagine that a proof of Conjecture \ref{conj1} would also apply to this much more general situation. See \cite{Ka2,Ka3} for more details on this conjecture.

One of the main problems with studying the effect of the minimal model programme on derived categories is the presence of singularities.
Considering varieties with mild singularities is essential for the minimal model programme in dimension at least three. But many of the techniques that have been developed to understand derived categories rely on smoothness.

The situation in dimension three is rather special since there is an explicit list of local models for terminal singularities.
Each is a finite quotient of a hypersurface singularity. The hypersurface singularity can  be thought of as a special fibre of a smooth fourfold,
which allows one to get around the bad behaviour of derived categories of singular varieties. In this way Chen \cite{Ch} was able to extend Theorem \ref{mine} so as to allow Gorenstein terminal singularities.

It is easy to see that if some version of Conjecture \ref{c3} is to hold, then for certain singular varieties the derived category must be modified in some way.
Thus for example there exist flops and flips which take one from a smooth variety to a singular one. 
But there can never be a fully faithful embedding $\D(Y_1)\into \D(Y_2)$
if $Y_1$ is singular and $Y_2$ smooth because of the following straightforward consequence of Serre's theorem on regularity of local rings.

\begin{prop}
Let $Y$ be a projective variety and  $\D(Y)$ its bounded derived category of coherent sheaves. Then
$Y$ is smooth if and only if
\[\dim_{\C}\bigoplus_{i\in\Z}\Hom^i_{\D(Y)}(A,B[i])<\infty\]
for all objects $A$ and $B$ of $\D(Y)$.
\end{prop}

At the moment this problem of how to modify the derived category for singular varieties is the biggest obstacle to progress in this area.

If a variety has only quotient singularities it is clear that one  should consider the derived category of coherent sheaves on the corresponding stack.
Kawamata explained this using the example of the Francia flip \cite{Ka1}. This is a threefold flip
\[
\xymatrix@C=.5em{ Y^+\ar[drr]_f&&&& Y^-\ar[dll]^g \\
&&X
}
\]
in which $Y^-$ is smooth but $Y^+$ has an isolated quotient singularity.
By the argument above there cannot be an embedding $\D(Y^+)\into\D(Y^-)$. On the other hand there is an embedding
$\D(\mathcal{Y}^+)\into\D(Y^-)$, where $\mathcal{Y}^+\to Y^+$ is the Deligne-Mumford stack associated to the quotient singularity.
In the same way, Chen's result on Gorenstein threefold flops extends to arbitrary terminal threefolds so long as one considers the derived categories of the corresponding stacks \cite{Ka2}.

Kawamata has  extended  Conjecture \ref{c3} to a statement involving log varieties with quotient singularities. In an impressive piece of work \cite{Ka3} he was able to prove his conjecture under the additional assumption that the birational maps $f_i$ were toroidal. Roughly speaking this means that they are locally defined by toric data. The derived McKay correspondence for abelian groups (Theorem \ref{abmc}) follows as a special case of this result.

Finally, Kawamata was able \cite{Ka4} to use his result  together with the log minimal model programme for toric varieties to solve a long-standing problem: the existence of a full exceptional collection  on any toric variety.  For more on these state of the art developments we refer the reader to \cite{Ka3,Ka4}.

\bigskip

\noindent Department of Pure Mathematics,
University of Sheffield,
Hicks Building, Hounsfield Road, Sheffield, S3 7RH, UK.

\smallskip

\noindent email: {\tt t.bridgeland@sheffield.ac.uk}

\end{document}